\documentclass[runningheads]{llncs}

\usepackage{amsmath, amsfonts, amssymb, graphicx, xcolor, dsfont, enumerate, dcolumn,subcaption}
\usepackage[numbers,square, comma, compress]{natbib}
\usepackage{algorithm}
\usepackage{algorithmic}
\usepackage{makecell}
\usepackage{threeparttable}

\let\doendproof\endproof
\renewcommand\endproof{\hfill\qed\doendproof}

\usepackage{hyperref}
\usepackage{xurl}
\usepackage{bm}
\usepackage[capitalise]{cleveref}
\PassOptionsToPackage{hyphens}
{url}\usepackage{hyperref}

\usepackage{tikz}

\newenvironment{theorem*}[1][Theorem]{\par\noindent\textbf{#1.}\itshape }{\par}
\newenvironment{corollary*}[1][Corollary]{\par\noindent\textbf{#1.}\itshape }{\par}

\begin{document}

\title{On the girth and connectivity of cubic graphs with a unique longest cycle}

\titlerunning{On the girth and connectivity of cubic graphs with a unique longest cycle}
%

\author{Jorik Jooken\inst{1} \and
Carol T. Zamfirescu\inst{2,3} }

\authorrunning{Jooken and Zamfirescu}

\institute{
Department of Computer Science, KU Leuven Kulak, 8500 Kortrijk, Belgium\\
\email{jorik.jooken@kuleuven.be}
\and
{Department of Applied Mathematics, Computer Science and Statistics, Ghent University, 9000 Ghent, Belgium}\\
\and
{
Department of Mathematics, Babe\c{s}-Bolyai University, Cluj-Napoca, Roumania
}
\email{czamfirescu@gmail.com}
}

\date{}

\maketitle

\begin{abstract}
We show that there exists an infinite family of cubic $2$-connected non-hamiltonian graphs with girth $5$ containing a unique longest cycle.

\end{abstract}


\section{Introduction}
The result proven in this note is fundamentally motivated by the problem of counting hamiltonian cycles in regular graphs, with a focus on hamiltonian regular graphs with few hamiltonian cycles. Perhaps the most important problem in this area is Sheehan's Conjecture~\cite{S75}, stating that every hamiltonian 4-regular graph has at least two distinct hamiltonian cycles. For recent articles treating the problem theoretically, see~\cite{BDP, GKN19, GJLSZ24, J24, Z22}, and for a recent algorithmic treatment, see~\cite{DMSZ24}. In this context we also want to highlight an interesting recent conjecture by Gir\~{a}o, Kittipassorn and Narayanan~\cite{GKN19} stating that there exists an absolute constant $K>0$ such that if an $n$-vertex graph $G$ with minimum degree at least $3$ contains a hamiltonian cycle, then $G$ contains another cycle of length at least $n-K$.

It has been known since the forties that if a 3-regular graph contains one hamiltonian cycle, then it must contain at least three~\cite{T46}. Chia and Thomassen explored in~\cite{CT12} a natural variation of this problem: instead of requiring the graph to contain a unique hamiltonian cycle, they investigated graphs with a unique \textit{longest cycle}. Although 3-regular graphs with a unique hamiltonian cycle do not exist, it is not difficult to describe 3-regular graphs containing exactly one longest cycle. For examples of such graphs satisfying various additional properties---e.g.\ regarding connectivity, girth, and length of the longest cycle---see for instance~\cite{CT12} itself and~\cite{Z22}.

Inspired by an old conjecture of Cantoni~\cite{T76} positing that every planar 3-regular graph with exactly three hamiltonian cycles contains a triangle, a question raised in~\cite{Z22} by the second author asks whether 3-regular 2-connected triangle-free graphs containing a unique longest cycle exist. Motivated, in turn, by this question, the first author posed a relaxation of this problem~\cite{J24} asking for the minimum number of vertices with a degree different from 3 amongst all 2-connected triangle-free non-hamiltonian graphs with minimum degree 3 containing a unique longest cycle. Moreover, he presented an example of a 3-regular 2-connected triangle-free non-hamiltonian graph containing precisely four longest cycles as well as a 2-connected triangle-free non-hamiltonian 60-vertex graph having a unique longest cycle and 56 vertices of degree 3 and 4 vertices of degree 4. In this note we answer in the affirmative the initial question, and thus also solve the proposed relaxation.

\section{The result}

\begin{theorem*}
\label{th:infFam}
If there exists a cubic $\kappa_1$-connected graph $G_1$ with girth $g_1$ having a unique longest cycle, a cubic graph $G_2$ with girth $g_2$ having precisely three hamiltonian cycles, and a vertex $u \in V(G_2)$ such that $G_2-u$ is $\kappa_2$-connected, then there exists an infinite family of cubic $\min\{\kappa_1,\kappa_2\}$-connected graphs with girth at least $\min\{2g_1,g_2\}$ having a unique longest cycle.
\end{theorem*}
\begin{proof}
Let $v, w$ and $x$ be the three neighbors of $u$ in $G_2$. For a graph $G$ and a vertex $a \in V(G)$ with precisely three neighbors $b$, $c$ and $d$, we define the \textit{marriage of vertex $a \in V(G)$ with $(G_2,u)$} as one of the graphs obtained by taking the disjoint union of $G-a$ and $G_2-u$ and adding three edges such that the resulting graph is cubic and every added edge has one endpoint in $\{b,c,d\}$ and the other endpoint in $\{v,w,x\}$ (for the purpose of this proof, it does not matter which of the six possible graphs we choose). Kotzig and Labelle~\cite{KL79} call this operation \textit{marriage}, whereas for example Eiben, Jajcay and Šparl~\cite{EJS19} refer to a very similar operation as \textit{generalized truncation}. For each vertex $t \in V(G_1)$, let $G_{2,t}$ be a copy of $G_2$ and let $u_t \in V(G_{2,t})$ be the vertex corresponding to $u \in V(G_2)$. Let $H$ be a graph obtained by marrying every vertex $t \in V(G_1$) with $(G_{2,t},u_t)$. We label the vertices of $H$ as $s_{i,j}$ in the natural way ($i \in V(G_1), j \in V(G_2-u)$).

$H$ is clearly cubic by construction. It is easy to check that for arbitrary vertices $p \in V(G_1)$ and $q \in V(G_2)$ if $G_1$ is $\kappa_1$-connected and $G_2-q$ is $\kappa_2$-connected, then the marriage of $p$ and $(G_2,q)$ is at least $\min\{\kappa_1,\kappa_2\}$-connected. Therefore $H$ is also $\min\{\kappa_1,\kappa_2\}$-connected. Let $\mathcal{C}=s_{i_0,j_0}s_{i_1,j_1} \ldots s_{i_{m-1},j_{m-1}}s_{i_0,j_0}$ be any cycle of length $m$ in $H$. If $i_0=i_1=\ldots=i_{m-1}$, then $j_0j_1 \ldots j_{m-1}j_0$ is a cycle of length $m$ in $G_2-u$. If $i_k \neq i_{k+1}$ for some $0 \leq k \leq m-1$, then $i_{k-1}=i_{k}$ and $i_{k+1}=i_{k+2}$ (where indices are taken modulo $m$), because $H$ is obtained by marrying the vertices of $G_1$ with $(G_2,u)$. Starting from the index $i_{k+1}$, generate a new sequence of indices $\mathcal{C}_2$ by moving sequentially through the indices, wrapping around from $i_{m-1}$ to $i_0$, until $i_{k+1}$ is reached again. At each step, include an index in $\mathcal{C}_2$ only if it is different from the immediately preceding index. Now $\mathcal{C}_2$ corresponds to a cycle in $G_1$. Hence, every cycle in $H$ corresponds to either a cycle in $G_2-u$ or a cycle in $G_1$ in which every vertex on that cycle is replaced by a path of length at least one. Clearly, the girth of $G_2-u$ is at least $g_2$. Therefore, the girth of $H$ is at least $\min\{2g_1,g_2\}$. 

Thomason~\cite{T78} showed that in any graph where all vertices have an odd degree every edge is contained in an even number of hamiltonian cycles. Because of Thomason's result, $G_2$ contains precisely one hamiltonian cycle containing simultaneously the edges $uv$ and $uw$, one hamiltonian cycle containing simultaneously the edges $uv$ and $ux$ and one hamiltonian cycle containing simultaneously the edges $uw$ and $ux$. As a result, the graph $G_2-u$ has precisely one hamiltonian $vw$-path, one hamiltonian $vx$-path and one hamiltonian $wx$-path. Hence, if $\mathcal{C}_1$ is the unique longest cycle in $G_1$, then the longest cycle in $H$ can have length at most $|\mathcal{C}_1| \cdot |V(G_2-u)|$ and there is indeed such a cycle (obtained by replacing each vertex in $\mathcal{C}_1$ by an appropriate hamiltonian path in $G_2-u$). Conversely, by contracting in $H$ every copy of $G_2-u$ to a single vertex, every longest cycle in $H$ corresponds to a longest cycle in $G_1$. Since $G_2-u$ has a unique hamiltonian $yz$-path for $y,z \in \{v,w,x\}, y \neq z$, we conclude that $H$ has a unique longest cycle.

By iteratively applying this argument, where $G_1$ is replaced by $H$, we obtain the desired infinite family.
\end{proof}

We remark that every cubic graph containing a unique longest cycle is necessarily non-hamiltonian due to Smith's observation in~\cite{T46}. In fact more strongly, one can infer that the length of the unique longest cycle in such a graph is at most $n-2$, where $n$ is the order of the graph, because Thomason~\cite{T78} proved that for every graph $G$ in which all vertices have an odd degree, the number of hamiltonian cycles in $G$ has the same parity as the number of hamiltonian cycles in $G-v$ for all $v \in V(G)$. 

\begin{figure}[ht]
\centering
\begin{tikzpicture}[scale=0.50]
\foreach \x in {0,1,...,8}{
\draw[fill] (\x*360/9:4.5) circle (1.8pt);
\draw[fill] (\x*360/9:3) circle (1.8pt);
\draw (\x*360/9:3)--(\x*360/9:4.5);
\draw (\x*360/9+360/9:4.5)--(\x*360/9:4.5);
\draw (\x*360/9+720/9:3)--(\x*360/9:3);
}
\end{tikzpicture}
\caption{A cubic 3-connected graph with girth 5 having precisely three hamiltonian cycles (this graph is known as the generalised Petersen graph GP(9,2)).}\label{fig:GMZ}
\end{figure}
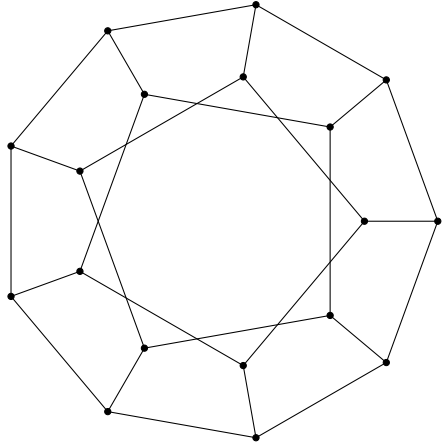

Chia and Thomassen~\cite{CT12} constructed a cubic 2-connected graph $G_1$ with girth 3 having a unique longest cycle (shown in Fig.~\ref{fig:ChiaAndThomassen}). Schwenk~\cite{S89} constructed a cubic 3-connected graph $G_2$ with girth 5 having precisely three hamiltonian cycles (shown in Fig.~\ref{fig:GMZ}). Therefore, for every vertex $u \in V(G_2)$, the graph $G_2-u$ is 2-connected. Note that if $G$ is a graph with girth at least $5$, then for any vertex $u \in V(G)$ the graph $G-u$ has girth at least $5$, but for our specific graph $G_2$, the graph $G_2-u$ has girth exactly $5$. Based on this observation and by applying the previous theorem, we obtain the following corollary.

\begin{figure}[ht]
\centering
\begin{tikzpicture}[scale=0.50]

\begin{scope}[rotate=179]
\foreach \x in {1,...,4}{
\draw[fill] (\x*360/5:3.25) circle (1.8pt);
\draw[fill] (\x*360/5-9:2.75) circle (1.8pt);
\draw[fill] (\x*360/5+9:2.75) circle (1.8pt);

\draw[fill] (\x*360/5:4.00) circle (1.8pt);
\draw[fill] (\x*360/5+5:4.5) circle (1.8pt);
\draw[fill] (\x*360/5-5:4.5) circle (1.8pt);

\draw (\x*360/5:3.25)--(\x*360/5-9:2.75);
\draw (\x*360/5:3.25)--(\x*360/5+9:2.75);
\draw (\x*360/5-9:2.75)--(\x*360/5+9:2.75);

\draw (\x*360/5:4.00)--(\x*360/5+5:4.5);
\draw (\x*360/5:4.00)--(\x*360/5-5:4.5);
\draw (\x*360/5-5:4.5)--(\x*360/5+5:4.5);

\draw (\x*360/5+720/5-9:2.75)--(\x*360/5+9:2.75);

\draw (\x*360/5:3.25)--(\x*360/5:4.00);
}

\draw (1*360/5+360/5-5:4.5)--(1*360/5+5:4.5);
\draw (3*360/5+360/5-5:4.5)--(3*360/5+5:4.5);
\draw (4*360/5+360/5:4.5)--(4*360/5+5:4.5);

\draw[fill] (0*360/5:3.25) circle (1.8pt);
\draw[fill] (0*360/5-9:2.75) circle (1.8pt);
\draw[fill] (0*360/5+9:2.75) circle (1.8pt);

\draw[fill] (0*360/5:4.5) circle (1.8pt);

\draw (0*60/5:3.25)--(0*360/5-9:2.75);
\draw (0*360/5:3.25)--(0*360/5+9:2.75);
\draw (0*360/5-9:2.75)--(0*360/5+9:2.75);

\draw (0*360/5+360/5-5:4.5)--(0*360/5:4.5);

\draw (0*360/5+720/5-9:2.75)--(0*360/5+9:2.75);

\draw (0*360/5:3.25)--(0*360/5:4.5);
\end{scope}

\begin{scope}[shift={(11,0)}, rotate=181+180]
\foreach \x in {1,...,4}{
\draw[fill] (\x*360/5:3.25) circle (1.8pt);
\draw[fill] (\x*360/5-9:2.75) circle (1.8pt);
\draw[fill] (\x*360/5+9:2.75) circle (1.8pt);

\draw[fill] (\x*360/5:4.00) circle (1.8pt);
\draw[fill] (\x*360/5+5:4.5) circle (1.8pt);
\draw[fill] (\x*360/5-5:4.5) circle (1.8pt);

\draw (\x*360/5:3.25)--(\x*360/5-9:2.75);
\draw (\x*360/5:3.25)--(\x*360/5+9:2.75);
\draw (\x*360/5-9:2.75)--(\x*360/5+9:2.75);

\draw (\x*360/5:4.00)--(\x*360/5+5:4.5);
\draw (\x*360/5:4.00)--(\x*360/5-5:4.5);
\draw (\x*360/5-5:4.5)--(\x*360/5+5:4.5);

\draw (\x*360/5+720/5-9:2.75)--(\x*360/5+9:2.75);

\draw (\x*360/5:3.25)--(\x*360/5:4.00);
}

\draw (1*360/5+360/5-5:4.5)--(1*360/5+5:4.5);
\draw (3*360/5+360/5-5:4.5)--(3*360/5+5:4.5);
\draw (4*360/5+360/5:4.5)--(4*360/5+5:4.5);

\draw[fill] (0*360/5:3.25) circle (1.8pt);
\draw[fill] (0*360/5-9:2.75) circle (1.8pt);
\draw[fill] (0*360/5+9:2.75) circle (1.8pt);

\draw[fill] (0*360/5:4.5) circle (1.8pt);

\draw (0*60/5:3.25)--(0*360/5-9:2.75);
\draw (0*360/5:3.25)--(0*360/5+9:2.75);
\draw (0*360/5-9:2.75)--(0*360/5+9:2.75);

\draw (0*360/5+360/5-5:4.5)--(0*360/5:4.5);

\draw (0*360/5+720/5-9:2.75)--(0*360/5+9:2.75);

\draw (0*360/5:3.25)--(0*360/5:4.5);

\end{scope}

\draw (7.15,-2.37)--(3.85,-2.37);
\draw (7.15,2.27)--(3.85,2.27);


\end{tikzpicture}
\caption{A cubic 2-connected graph with girth 3 having a unique longest cycle.}
\label{fig:ChiaAndThomassen}
\end{figure}
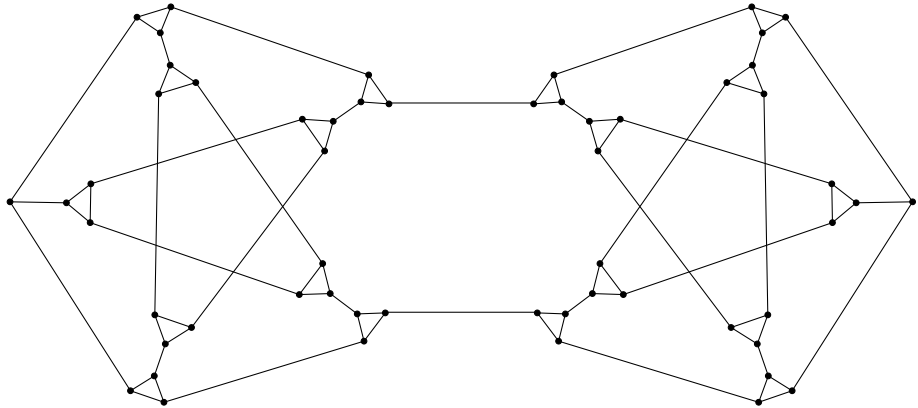

\vspace{0.3cm}
\begin{corollary*}
\label{cor:result}
There exists an infinite family of cubic $2$-connected non-hamiltonian graphs with girth $5$ containing a unique longest cycle.
\end{corollary*}


\section{Two open problems}

We conclude this note with two questions that spark our interest.

\begin{question}
Does there exist a cubic $2$-connected non-hamiltonian \textit{planar} triangle-free graph containing a unique longest cycle?
\end{question}
If Cantoni's conjecture~\cite{T76} holds, it seems that a technique different from the one presented in the theorem will be necessary to construct such graphs if they exist.

The second question that interests us is as follows:

\begin{question}
    Does there exist a cubic $3$-connected graph with girth at least $g$ having precisely three hamiltonian cycles for each integer $g$?
\end{question}
If the answer is positive, one can infer the existence of an infinite family of cubic 2-connected non-hamiltonian graphs with girth at least $g$ containing a unique longest cycle for each integer $g$ (by using the corollary and then iteratively applying the theorem).

\section*{Acknowledgements}
Jorik Jooken is supported by a Postdoctoral Fellowship of the Research Foundation Flanders (FWO) with grant number 1222524N.

\end{document}